\RequirePackage{ifpdf}
\ifpdf % We are running pdfTeX in pdf mode
\documentclass[pdftex]{sigma}
\else
\documentclass{sigma}
\fi

\begin{document}

\allowdisplaybreaks

\renewcommand{\PaperNumber}{085}

\FirstPageHeading

\ShortArticleName{On the Projective Algebra of Randers Metrics of Constant Flag Curvature}

\ArticleName{On the Projective Algebra of Randers Metrics\\ of Constant Flag Curvature}

\Author{Mehdi RAFIE-RAD~$^{\dag\ddag}$ and Bahman REZAEI~$^\S$}

\AuthorNameForHeading{M.~Rafie-Rad and B.~Rezaei}

\Address{$^\dag$~School of Mathematics, Institute for Research in Fundamental Sciences (IPM),\\
\hphantom{$^\dag$}~P.O. Box 19395-5746, Tehran, Iran}

\Address{$^\ddag$~Department of Mathematics, Faculty of Mathematical Sciences, University of Mazandaran,\\
\hphantom{$^\ddag$}~P.O. Box 47416-1467, Babolsar, Iran}
\EmailD{\href{mailto:rafie-rad@umz.ac.ir}{rafie-rad@umz.ac.ir}, \href{mailto:m.rafiei.rad@gmail.com}{m.rafiei.rad@gmail.com}}

\Address{$^\S$~Department of Mathematics, Faculty of Sciences,
University of Urmia, Urmia, Iran}
\EmailD{\href{mailto:b.rezaei@urmia.ac.ir}{b.rezaei@urmia.ac.ir}}

\ArticleDates{Received February 26, 2011, in f\/inal form August 20, 2011;  Published online August 31, 2011}

\Abstract{The collection of all projective vector f\/ields on a Finsler space $(M, F)$ is a f\/inite-dimensional Lie algebra with respect to the usual Lie bracket, called the projective algebra denoted by $p(M,F)$ and is the Lie algebra of the projective group $P(M,F)$. The projective algebra
$p(M,F=\alpha+\beta)$ of a Randers space is characterized as a certain Lie subalgebra of the projective algebra $p(M,\alpha)$. Certain subgroups of
the projective group $P(M,F)$ and their invariants are studied. The projective algebra of Randers metrics of constant f\/lag curvature is studied and it is proved that the dimension of the projective algebra of Randers metrics constant f\/lag curvature on a compact $n$-manifold either equals $n(n+2)$ or at most is $\frac{n(n+1)}{2}$.}

\Keywords{Randers metric; constant f\/lag curvature; projective vector f\/ield; projective al\-gebra}

\Classification{53C60; 53B50, 58J60}

\section{Introduction}

The motion of freely falling particles def\/ine a projective structure on spacetime. Mathematically speaking, this provides a projective connection or an equivalence class of symmetric af\/f\/ine connections all possessing the same unparameterized geodesic curves. This may be regarded as a~mathematical formulation of the weak principle of equivalence valid both in the Newtonian and relativistic theory of spacetime and gravity~\cite{Israel1973}. Many physical considerations require metric structures on spacetime in liaison to af\/f\/ine connections. A necessary condition for two such metric structures to have the same (unparameterized) geodesic curves is that their Weyl projective tensors are identical.

The locally anisotropic space-times are studied in \cite{Stavrinos2009} from a geometrical point of view and thus may include some auspices on the Weyl projective tensor. As we will see in this paper, certain subgroups of the Lorentz group may be at once a subgroup of the projective group of the Finsler metric $F=\sqrt{\eta_{\mu\nu}dx^\mu dx^\nu}+{\bf A}_\mu dx^\mu$, where $\eta$ and ${\bf A}$ denotes the Lorentz metric and the electromagnetic potential vector of the f\/lat space-time. Some reduced forms of Weyl projective tensor $W$ have been introduced in \cite{Akbar-Zadeh1986,NajafiTayebi2010}, which are invariant among projectively related constant curvature Finsler metrics but not identical among scalar f\/lag curvature metrics. Any two Finsler metrics possessing the same (unparameterized) geodesics have the same Weyl projective tensor. Studying Weyl projective vector f\/ields (i.e.\ those vector f\/ields preserving the Weyl projective tensor and also the reduced Weyl projective tensors) and projective vector f\/ields have such a leading role to obtain projective symmetries which provide some conservation laws in physical terms. On the other hand, there are many papers devoted on projective symmetry in metric-af\/f\/ine gravity and cosmology, see for example \cite{HallLonie2008,Barnes1993}.

Randers metrics are the most popular Finsler metrics in dif\/ferential geometry and physics simply obtained by a Riemannian metric $\alpha=\sqrt{a_{ij}(x)y^iy^j}$ and $\beta=b_i(x)y^i$ as $F=\alpha+\beta$ and was introduced by G.~Randers in \cite{Randers1941} in the context of general relativity. They arise naturally as the geometry of light rays in stationary spacetimes \cite{GibbonsHerdeiroWerner2009}. One may refer to \cite{BaoRobles2003,BaoShen2002,Mo2008,Robles2003} for an extensive series of results about the Einstein Randers metrics and the Randers metrics of constant f\/lag curvature. The present paper is closely related to the problem of projective relatedness of Randers metrics which is investigated in~\cite{ShenYu2008}. To avoid the obscurity, given a~Randers metric $\alpha+\beta$, the geometric objects in $(M,F)$ and $(M,\alpha)$ are denoted respectively by the prif\/ices ``$F$-'' and ``$\alpha$-'', respectively, for instance, an $F$-projective vector f\/ield means a~projective vector f\/ield on $(M,F)$, an $\alpha$-projective vector f\/ield means a~projective vector f\/ield on $(M,\alpha)$, an $\alpha$-Killing vector f\/ields stands for a Killing vector f\/ield on $(M,\alpha)$, etc. We use the usual notations for Randers metrics in \cite{BaoShen2002,ShenYu2008}. Given any vector f\/ield $V$, its complete lift to $TM_0=TM\backslash\{0\}$ is denoted by $\hat{V}$ and the Lie derivative along $\hat{V}$ is denoted by ${\cal L}_{\hat{V}}$. One can benef\/it~\cite{Yano1957} for an extensive discussion on the theory of Lie derivatives of various geometric objects in Finsler spaces. Traditionally we use the notations for the so-called $(\alpha,\beta)$-metric disposal in~\cite{Shen2001}, thereby $s^i_{\ \circ}=a^{ik}s_{kj}y^j$, where $s_{kj}=\big(\frac{\partial b_k}{\partial x^j}-\frac{\partial b_j}{\partial x^k}\big)$. We characterize the projective vector f\/ields on a Randers space by proving the following theorem:

\begin{theorem}
\label{Liethm}
Let $(M,F=\alpha+\beta)$ be a Randers space and $V$ be a vector field on $M$. $V$ is $F$-projective if and only if $V$ is $\alpha$-projective and
${\cal L}_{\hat{V}}\{\alpha s^i_{\ \circ}\}=0$.
\end{theorem}

Determining the dimension of the projective algebra of constant curvature and Einstein spaces is of interests in physical and geometrical discussions,
see \cite{Barnes1993} and interested readers may be\-nef\/it~\cite{Yano1957} for a large discussion on this f\/ield. This leads to calculate number of independent projective vector f\/ields and is closely related to the number of independent Killing vector f\/ields in each case. It is well known that in an $n$-dimensional Riemannian space of constant curvature the dimension of the projective algebra is $n(n+ 2)$ and vice-versa, see \cite{Barnes1993,Yano1957}. This weaves an overture for an analogue problem for Randers space $(M,F=\alpha+\beta)$. If we have $s^i_{\ j}=0$, then the  respective projective algebras $p(M,F)$ of $F$ and $p(M,\alpha)$ of $\alpha$ coincide. If moreover, $F$ is locally projectively f\/lat, then $\alpha$ is too and hence, the dimension of the projective algebra $p(M,F)$ is $n(n+2)$. Notice
that our discussions are closely related to the algebra $k(M,\alpha)$ of $\alpha$-Killing vector f\/ields. The important case is considerable when
$s^i_{\ j}\neq0$ and uncovers a non-Riemannian feature of Finsler metrics in comparison with the analogue Riemannian case. We summarize the argument by establishing the following result:

\begin{theorem}
\label{admitting}
Let $(M,F)$ be an $n$-dimensional $(n\geq3)$ Randers space of constant flag curvature and $M$ is compact. The dimension of the projective algebra $p(M,F)$ is either $n(n+2)$ or at most equals $\frac{n(n+1)}{2}$.
\end{theorem}

The spaces admitting certain vector f\/ields has a long history in Riemannian geometry,  see for example \cite{Akbar-Zadeh1978,Barnes1993,Hiramoto1980,Obata1962,Tanno1978,Tashiro1965,Yamauchi1974}. Existence of some special vector f\/ields on a
Riemannian space may pertain to some global properties of the underlying Riemannian space. We prove the following result to uncover such an interaction for Randers spaces:
\begin{theorem}
\label{mainthm}
Let $(M,F=\alpha+\beta)$ be a Randers space of vanishing ${\bf S}$-curvature and dimension $n\geq3$. If $(M,F)$ admits a non-$\alpha$-affine projective vector field $V$, then
$(M,F)$ is a Berwald space.
\end{theorem}

\section{Preliminaries}\label{sectionP}

Let $M$ be a $n$-dimensional $ C^\infty$ connected manifold.  $T_x M $ denotes  the tangent space of~$M$ at~$x$. The tangent bundle of $M$ is the union of tangent spaces $TM:=\bigcup _{x \in M} T_x M$.  We will denote  the elements of $TM$ by $(x, y)$ where $y\in T_xM$. Let $TM_0 = TM\setminus \{ 0 \}.$ The natural projection $\pi: TM_0 \rightarrow M$ is given by $\pi (x,y):= x$.  A~\textit{Finsler metric} on  $M$ is a function $ F:TM \rightarrow [0,\infty )$ with the following properties: $(i$)~$F$~is $C^\infty$ on $TM_0$,  $(ii)$~$F$~is positively 1-homogeneous on the f\/ibers of tangent bundle $TM$,  and  $(iii)$~the Hessian of $F^{2}$ with elements $ g_{ij}(x,y):=\frac{1}{2}[F^2(x,y)]_{y^iy^j} $ is positive def\/inite matrix on $TM_0$.  The pair  $(M,F)$ is then called a {\it Finsler space}. Throughout this paper, we denote a Riemannian metric by  $\alpha=\sqrt{a_{ij}(x)y^iy^j}$ and a 1-form by $\beta=b_i(x)y^i$.
A globally def\/ined vector f\/ield ${\bf G}$ is induced by $F$ on $TM_0$, which in a standard coordinate $(x^i,y^i)$
for $TM_0$ is given by ${\bf G}=y^i {{\partial} \over {\partial x^i}}-2G^i(x,y){{\partial} \over {\partial y^i}}$, where $G^i(x,y)$ are local functions on $TM_0$ satisfying $G^i(x,\lambda y)=\lambda^2 G^i(x,y)$, $\lambda>0$. Assume the following conventions:
\[
G^i_{\ j}=\frac{\partial G^i}{\partial y^i},\qquad
G^i_{\ jk}=\frac{\partial G^i_{\ j}}{\partial y^k},\qquad
G^i_{\ jkl}=\frac{\partial G^i_{\ jk}}{\partial y^l}.
\]
Notice that the local functions $G^i_{\ jk}$ give rise to a torsion-free connection in $\pi^*TM$ called \textit{the Berwald connection} which is practical in this paper, see~\cite{Shen2001}. The local functions $G^i_{\ j}$ def\/ine a~nonlinear connection ${\cal H}TM$ spanned by the horizontal frame $\{\frac{\delta}{\delta x^i}\}$, where $\frac{\delta}{\delta x^j}=\frac{\partial}{\partial x^j}-G^i_{\ j}\frac{\partial}{\partial y^i}$. The nonlinear connection ${\cal H}TM$ splits $TTM$ as $TTM=\ker\pi_*\oplus{\cal H}TM$, see~\cite{Shen2001}. A Finsler metric is called a {\it Berwald metric} if $G^i_{\ jk}(x, y)$ are functions of $x$ only at every point $x\in M$, equivalently $F$ is a Berwald metric if and only if $G^i_{\ jkl}=0$.

For a Finsler metric $F$ on an $n$-dimensional manifold $M$ {\it the
Busemann--Hausdorff volume form} $dV_F = \sigma_F(x) dx^1 \cdots
dx^n$ is def\/ined by
\[
\sigma_F(x) := \frac{\textrm{Vol} ({\mathbb B}^n(1))}{
\textrm{Vol}  \{ (y^i)\in \mathbb{R}^n \;  | \; F  ( y^i
\frac{\partial}{\partial x^i}|_x  ) < 1  \} }.
\]
Assume $\underline{g}={\det  ( g_{ij}(x,y)  )}$ and def\/ine $\tau (x, y):=\ln{\sqrt{\underline{g}}\over \sigma_F(x)}$. $\tau=\tau(x,y)$ is a scalar function on~$TM_0$, which is called the {\it distortion}~\cite{Shen2001}. For a vector $y\in T_xM$, let $c(t)$, $-\epsilon < t <\epsilon $, denote the
geodesic with $c(0)=x$ and $\dot{c}(0)=y$. The function ${\bf S}(y):= {d \over dt}  [ \tau  (\dot{c}(t) )  ]_{|_{t=0}}$ is called the ${\bf S}$-curvature with respect to Busemann--Hausdorf\/f volume form. A Finsler space is said to be {\it of isotropic ${\bf S}$-curvature} if there is a function $\sigma=\sigma(x)$ def\/ined on $M$ such that ${\bf S}=(n+1)\sigma(x)F$. It is called a Finsler space {\it of constant {\bf S}-curvature} once $\sigma$ is a constant. Every Berwald space is of vanishing {\bf S}-curvature \cite{Shen2001}. The {\bf E}-curvature of the Finsler space $(M,F)$ is def\/ined by ${\bf E}_y={\bf E}_{ij}(y)dx^i\otimes dx^j$, where ${\bf E}_{ij}=\frac{1}{2}\frac{\partial^2{\bf S}}{\partial y^i\partial y^j}$. $(M,F)$ is called a weakly-Berwald space if ${\bf E}=0$. It is easy to see that we have ${\bf E}_{ij}=\frac{1}{2}G^r_{\ irj}$.

Let $(M,\alpha)$ be a Riemannian space and $\beta=b_i(x)y^i$ be a 1-form def\/ined on $M$ such that $\|\beta\|_x :=\sup\limits_{y \in T_xM} \beta(y)/\alpha(y) < 1$. The Finsler metric $F = \alpha+\beta$ is called a Randers metric on a manifold $M$. Denote the geodesic spray coef\/f\/icients of $\alpha$ and $F$ by
the notions  $G_\alpha^i$ and $G^i$, respectively and the Levi-Civita connection of $\alpha$ by $\nabla$. Def\/ine $\nabla_jb_i$ by $(\nabla_jb_i) \theta^j := db_i -b_j \theta_i^{\ j}$, where $\theta^i :=dx^i$ and $\theta_i^{\ j} :=\tilde{\Gamma}^j_{ik} dx^k$ denote the Levi-Civita  connection forms and $\nabla$ denotes its associated covariant derivation of $\alpha$. Let us put
\begin{gather*}
r_{ij} := {1\over 2}  ( \nabla_jb_i+\nabla_ib_j), \qquad s_{ij}:= {1\over 2}(\nabla_jb_i-\nabla_ib_j),\\
  s^i_{\ j} :=  a^{ih}s_{hj}, \qquad s_j:=b_i s^i_{\ j}, \qquad e_{ij} := r_{ij}+ b_i s_j + b_j s_i.
\end{gather*}
Then $G^i$  are given by
\begin{gather*}
G^i = G_\alpha^i + \left({e_{\circ\circ} \over 2F} -s_\circ\right)y^i+ \alpha s^i_{\ \circ},%\label{Gi}
\end{gather*}
where $e_{\circ\circ}:= e_{ij}y^iy^j$, $s_\circ:=s_iy^i$, $s^i_{\ \circ}:=s^i_{\ j} y^j$
and $G_\alpha^i$ denote the geodesic coef\/f\/icients of $\alpha$,
see~\cite{Shen2001}.
Notice that the {\bf S}-curvature of a Randers metric $F=\alpha+\beta$ can be obtained as follows
\begin{gather*}
%\label{S formula}
{\bf S}=(n+1)\left\{\frac{e_{\circ\circ}}{F}-s_\circ-\rho_\circ\right\},
\end{gather*}
where $\rho=\ln \sqrt{1-\|\beta\|}$ and $\rho_\circ=\frac{\partial \rho}{\partial x^k}y^k$. It is well-known that every weakly-Berwald Randers space is of vanishing {\bf S}-curvature~\cite{Shen2001}.

Let $F$ be a Finsler metric on an $n$-manifold and $G^i$ denote the geodesic coef\/f\/icients of $F$. Def\/ine
${\bf R}_y= K^i_{\ k}(x, y) dx^k \otimes \frac{\partial}{\partial x^i}|_x: T_xM \to T_xM$ by
\[
 K^i_{\ k} := 2 \frac{\partial G^i}{\partial x^k}
- y^j \frac{\partial^2 G^i}{\partial x^j\partial y^k} + 2 G^j \frac{\partial^2 G^i}{\partial y^j \partial y^k}
- \frac{\partial G^i}{\partial y^j} \frac{\partial G^j}{\partial y^k}.
\]
The family ${\bf R}:=\{{\bf R}_y\}_{y\in TM_0}$ is called the Riemann curvature~\cite{Shen2001}. The {\it Ricci scalar} is denoted by ${\bf Ric}$ it is def\/ined by ${\bf Ric}:=K^k_{\ k}$. The Ricci scalar ${\bf Ric}$ is a generalization of the Ricci tensor in Riemannian geometry. A~Finsler space $(M,F)$ is called an {\it Einstein space} if there is function~$\sigma$ def\/ined on $M$ such that ${\bf Ric}=\sigma(x)F^2$. D.~Bao and C.~Robles proved in~\cite{BaoRobles2003,Robles2003} the
following theorem:

\begin{theorem}   Let $(M,F=\alpha+\beta)$ be an $n$-dimensional Randers space and $n\geq3$. If $(M,F)$ is an Einstein space
with ${\bf Ric} = (n-1)K(x)F^2$, then it is of constant ${\bf S}$-curvature and $K(x)$ is constant.
\end{theorem}

The Berwald--Riemannian curvature tensor ${\bf K}_y=K^i_{\ jkl}(y)\frac{\partial}{\partial x^i}\otimes dx^j\otimes dx^k\otimes dx^l$ and the Berwald--Ricci tensor $K_{jl}(y)dx^j\otimes dx^l$ are respectively def\/ined  by
\begin{gather*}
K^i_{\ jkl} := \frac{1}{3}\left\{\frac{\partial^2K^i_{\ k}}{\partial y^j\partial y^l}-\frac{\partial^2K^i_{\ l}}{\partial y^j\partial y^k}\right\},\qquad K_{jl}:=K^i_{\ jil}.
\end{gather*}
Due to a result in \cite{Mo2008}, every Finsler metric of constant {\bf S}-curvature on a compact manifold is of vanishing {\bf S}-curvature. Therefore, the {\bf S}-curvature of every Einstein Randers metric on an $n$-dimensional ($n\geq3$) compact manifold is vanishing.

Denote the horizontal and vertical covariant derivation of the Berwald connection of $F$ respectively by ``$_|$" and ``$_.$". The quit nouvelle non-Riemannian quantity ${\bf H}_y={\bf H}_{ij}(y)dx^i\otimes dx^j$ is simply def\/ined by ${\bf H}_{ij}={\bf E}_{ij|k}y^k$, see \cite{Akbar-Zadeh1988,Mo2009,NajafiTayebiShen2008}. Consider the following \textit{Bianchi identity} for the Berwald connection \cite{Akbar-Zadeh1988}:
\[
G^i_{\ jkl|m}-G^i_{\ jkm|l}=K^i_{\ jkl.m}.\label{Ricci identity}
\]
After convecting the indices $i$ and $k$ and taking into account the equation $G^i_{\ jil}=2{\bf E}_{jl}$, we obtain $G^k_{\ jkl|m}-G^k_{\ jkm|l}=2({\bf E}_{\ jl|m}-{\bf E}_{\ jm|l})=K_{\ jl.m}$. From which it results
\begin{gather}
y^jK_{\ jl.m}=0,\qquad  y^lK_{\ jl.m}=-2{\bf H}_{jm}.\label{Ricci identity}
\end{gather}

\section{Projectively related metrics and projectively invariants}

Two Finsler metrics $F$ and $\tilde{F}$ on a manifold $M$ are said to be \textit{$($pointwise$)$ projectively related} if
they have the same geodesics as point sets. Hereby, there is a function $P(x,y)$ def\/ined on $TM_0$ such that $\tilde{G}^i=G^i+Py^i$ on
coordinates $(x^i,y^i)$ on $TM_0$, where $\tilde{G}^i$ and $G^i$ are the geodesic spray coef\/f\/icients of $\tilde{F}$ and $F$, respectively. A Finsler metric $F$ on an open subset $\textsl{U}\subseteq\mathbb{R}^n$ is called \textit{projectively flat} if all geodesics are straight in~$\textsl{U}$. In this case, $F$ and the Euclidean metric on $\textsl{U}$ are projectively related. A Finsler metric is called \textit{locally projectively} f\/lat if at any point $x\in M$, there is a local coordinate $(x^i, U)$ in which $F$ is projectively f\/lat. We consider projectively related Finsler metrics, namely those having the same geodesics as set points. Let~$\tilde{F}$ and~$F$ be two projectively related Finsler metrics. Consider a natural coordinate system $((x^i,y^i),\pi^{-1}(U))$. There is function~$P$ def\/ined on $TM_0$ such that $\widetilde{G}^i=G^i+Py^i$. Let us put $P_i=P_{.i}$ and $P_{ij}=P_{i.j}$. Observe that we have
\begin{gather}
\widetilde{G}^i_{\ j} = G^i_{\ j}+P_jy^i+P\delta^i_{\ j},\qquad \widetilde{G}^i_{\ jk}=G^i_{\ jk}+P_{jk}y^i+P_{k}\delta^i_{\ j}+P_{j}\delta^i_{\
k},
\label{1}\\
\widetilde{\bf E}_{ij} = {\bf E}_{ij}+\frac{(n+1)}{2}P_{ij}.\label{3}
\end{gather}
The Berwald--Riemannian curvature and the Berwald--Ricci tensors of $\tilde{F}$ and $F$ are related as follows
\begin{gather}
\widetilde{K}^i_{\,\,hjk} =  K^i_{\,\,hjk}+y^i(P_{jh|k}-P_{kh|j})+\delta^i_h(P_{j|k}-P_{k|j})\nonumber\\
\phantom{\widetilde{K}^i_{\,\,hjk} =}{} +\delta^i_j(P_{h|k}-P_hP_k-PP_{hk})-\delta^i_k(P_{h|j}-P_hP_j-PP_{hj}),\nonumber\\
\widetilde{K}_{\,\,hk} = K_{hk}+(P_{h|k}-P_{k|h})+(n-1)(P_{h|k}-P_hP_k-PP_{hk})-P_{hk|\circ}.\label{RicK}
\end{gather}
Finally  we f\/ind out that $(\widetilde{K}_{\,\,hk}-\widetilde{K}_{\,\,kh})_{.j}=(K_{hk}-K_{kh})_{.j}+(n+1)(P_{h|k}-P_{k|h})_{.j}$.
The non-Riemannian quantity {\bf H} was introduced in~\cite{Akbar-Zadeh1986,NajafiTayebiShen2008} and developed in~\cite{Mo2009}. We would like to consider projectively related Finsler metrics with the same {\bf E}- and {\bf H}-curvatures. Observe that according to~(\ref{1}) and~(\ref{3}), {\bf
H}-curvatures of $\tilde{F}$ and $F$ are related as follows
\begin{gather}
\widetilde{{\bf H}}_{ij}=
y^r\frac{\tilde{\delta}}{\tilde{\delta}x^r}\widetilde{\bf E}_{ij}-\widetilde{\bf E}_{rj}\widetilde{G}^r_{\ i}-\widetilde{\bf
E}_{ir}\widetilde{G}^r_{\ j}=\left(y^r\frac{\delta}{\delta x^r}-2Py^r\frac{\partial}{\partial y^r}\right)
\left({\bf E}_{ij}+\frac{(n+1)}{2}P_{ij}\right)\nonumber\\
\phantom{\widetilde{{\bf H}}_{ij}}{} -\left({\bf E}_{rj}+\frac{(n+1)}{2}P_{rj}\right)\!(G^r_{\ i}+P_iy^r+P\delta^r_{\ i})-
\left({\bf E}_{ir}+\frac{(n+1)}{2}P_{ir}\right)\!(G^r_{\ j}+P_jy^r+P\delta^r_{\ j})\nonumber \\
\phantom{\widetilde{{\bf H}}_{ij}}{} = {\bf E}_{ij|\circ}+\frac{(n+1)}{2}P_{ij|\circ}
={\bf H}_{ij}+\frac{(n+1)}{2}P_{ij|\circ}, \label{Hc}
\end{gather}
where $\frac{\tilde{\delta}}{\tilde{\delta} x^k}=\frac{\partial}{\partial x^k}-2\widetilde{G}^i_{\ k}\frac{\partial}{\partial y^i}$ and
$\frac{\delta}{\delta x^k}=\frac{\partial}{\partial x^k}-2G^i_{\ k}\frac{\partial}{\partial y^i}$.
From (\ref{Hc}) one may conclude the following lemma.

\begin{lemma}
\label{Lemma 1}
Suppose that $\widetilde{F}$ and $F$ are projectively related with projective factor $P$. Then~$\widetilde{F}$ and~$F$ have the same {\bf H}-curvature if and only if $P_{ij|\circ}=0$.
\end{lemma}

The Finsler metrics $\widetilde{F}$ and $F$ are said to be \textit{{\bf H}-projectively related} if they are projectively related and have the same {\bf H}-curvature. From (\ref{3}) it results that if given any $x\in M$ the function~$P(x,y)$ is linear with respect to $y$, in other words $P_{ij}=0$, then $F$ and $\widetilde{F}$ are {\bf H}-projectively related. Hereby the Finsler metrics $\widetilde{F}$ and $F$ are said to be \textit{specially projectively related} if~$P(x,y)$ is linear with respect to $y$.

\begin{example}
The Funk metric $\Theta$ on the Euclidean unit ball $\mathbb{B}^n(1)$ is a Randers metric given by
\[
  \Theta(x,y):=\frac{\sqrt{|y|^2-(|x|^2|y|^2-\langle x,y\rangle^2)}}{1-|x|^2}+\frac{\langle x,y\rangle}{1-|x|^2}, \qquad y\in T_xB^n(1)\simeq \mathbb{R}^n,
\]
where $\langle,\rangle$ and $|.|$ denotes the Euclidean inner product and norm on $\mathbb{R}^n$, respectively. Given any constant vector $a\in \mathbb{R}^n$, the generalized Funk metric $\Theta_a$ is given by $\Theta_a:=\Theta+d\varphi_a$, where $\varphi_a=\ln (1+\langle a,x\rangle)+C$ and $C$ is a constant. From the variational point of view this changes the length function by something which depends only on the end-points, not the path between them. One may also refer to \cite{Shen2001} to f\/ind an analytic proof. $\Theta$ and $\Theta_a$ are both projectively f\/lat, {\bf H}-projectively related of constant {\bf S}-curvature with $\sigma=\frac{1}{2}$. It not hard to see that the projective factor $P$ is $P=-\frac{1}{2}d\varphi_a$.
\end{example}

\begin{example}
Given any vector $a\in\mathbb{R}^n$, def\/ine the Finsler metric $F$ on $\mathbb{B}^n(1)$ by
\[
F := (1+\langle a,x\rangle)\big(\Theta+\Theta_{x^k}x^k\big).
\]
$F$ is projectively f\/lat with projective factor $P =\Theta$ and $F$ is of constant f\/lag curvature ${\bf K} = 0$. Thus $F$ and the Euclidean metric on
$\mathbb{B}^n(1)$ have the same vanishing {\bf H}-curvature and are {\bf H}-projectively related. This example is borrowed from~\cite{Shen2002}.
\end{example}

\subsection{Projectively invariants}

Any geometric object which is identical between two projectively related metrics is called a~\textit{projective invariant}. There are many
projectively invariant tensors in Finsler geometry such as the {\it Douglas} tensor ${\bf D}_y=D^i_{\ jkl}(y)\frac{\partial}{\partial x^i}\otimes dx^j\otimes dx^k\otimes dx^l$ and \textit{Weyl} curvature ${\bf W}_y=W^i_{\ jkl}(y)\frac{\partial}{\partial x^i}\otimes dx^j\otimes dx^k\otimes dx^l$. However, the notion of the projective connection in Finsler geometry encounters some dif\/f\/iculties to be globally projectively invariant. The tensors ${\bf D}$ and ${\bf W}$ are def\/ined as follows
\begin{gather*}
D^i_{\ jkl}=\frac{\partial^3}{\partial y^j\partial y^k\partial y^l}\left\{G^i-\frac{1}{n+1}G^m_{\ m}y^i\right\},\nonumber\\
W^i_{\ jkl}=K^i_{\ jkl}-\frac{1}{n^2-1}\big\{\delta^i_{\
j}(\hat{K}_{kl}-\hat{K}_{lk})+(\delta^i_{\ k}\hat{K}_{jl}-\delta^i_{\
l}\hat{K}_{jk})+y^i(\hat{K}_{kl}-\hat{K}_{lk})_{.j}\big\},%\label{Douglas}
\end{gather*}
where $\hat{K}_{jk}=nK_{jk}+K_{kj}+y^rK_{kr.j}$. In 1986, H.~Akbar-Zadeh introduced a tensor which is just invariant by a~sub-group of projective
transformations, not all of them~\cite{Akbar-Zadeh1988}. In fact, this is a~non-Riemannian generalization of Weyl's curvature. It is denoted by
$\overset{*}{W}{^i_{\ jkl}}$ and is def\/ined by
\[
\overset{*}{W}{^i_{\ jkl}}=K^i_{\ jkl}-\frac{1}{n^2-1}\big\{\delta^i_{\
k}(nK_{jl}+K_{lj})-\delta^i_{\ l}(nK_{jk}+K_{kj})+(n-1)\delta^i_{\ j}(K_{kl}-K_{lk})\big\}.
\]
Assume that $W^i_{\ k}=y^jy^lW^i_{\ jkl}$. From (\ref{Ricci identity}) ${\bf W}$ can be written in terms of {\bf H}-curvature by the following equation
\begin{gather}
\label{W H}
W^i_{\ jkl}=\overset{*}{W}{^i_{\ jkl}}-\frac{2}{n^2-1}\big\{\delta^i_{\ l}{\bf H}_{jk}-\delta^i_{\ k}{\bf H}_{jl}\big\}
-\frac{y^i}{n+1}(K_{kl}-K_{lk})_{.j}.
\end{gather}
One may easily check from (\ref{RicK}) and (\ref{Hc}) that every two specially projectively related metrics $\tilde{F}$ and $F$ have the same tensors $W$, {\bf H} and $(K_{hk}-K_{kh})_{.j}$. Observe that from (\ref{W H}) it results that they have the same tensor $\overset{*}{W}$. There is the following identity for $W$ given in~\cite{Akbar-Zadeh1986,Shen2001}:
\begin{gather}
\label{W identity}
W^i_{\ jkl}=\frac{1}{3}\big\{W^i_{\ k.l.j}-W^i_{\ l.k.j}\big\}.
\end{gather}

%\noindent
%{\bf Theorem B.}
\begin{theorem}
Let $(M,F)$ be an $n$-dimensional Finsler manifold $(n\geq3)$. ${\bf W}=0$ if and only if~$F$ is of scalar flag curvature.
\end{theorem}

Let $\overset{*}{W}{^i_{\ k}}=y^jy^l\overset{*}{W}{^i_{\ jkl}}$. The following theorem is proved in \cite{Akbar-Zadeh1986}, however, we give a modif\/ied proof for
it.

%\noindent
%{\bf Theorem C.}
\begin{theorem}
Let $(M,F)$ be an $n$-dimensional Finsler manifold $(n\geq3)$. $\overset{*}{W}=0$ if and only if~$F$ is of constant flag
curvature.
\end{theorem}

\begin{proof}
From (\ref{W identity}), (\ref{W H}) and $y^l{\bf H}_{jl}=0$, it follows that we have $W^i_{\ k}=\overset{*}{W}{^i_{\ k}}$ and
\[
W^i_{\ jkl}=\frac{1}{3}\big\{\overset{*}{W}{^i_{\ k.l.j}}-\overset{*}{W}{^i_{\ l.k.j}}\big\}.
\]
Now let $\overset{*}{W}=0$. It follows immediately that $W=0$ and from Theorem B. that $(M,F)$ is of scalar curvature. But, from (\ref{W H}) it
results
\[
\frac{2}{n^2-1}\big\{\delta^i_{\ l}{\bf H}_{jk}-\delta^i_{\ k}{\bf H}_{jl}\big\}+\frac{y^i}{n+1}(K_{kl}-K_{lk})_{.j}=0.
\]
Convecting the index $k$ in $y^k$ and applying (\ref{Ricci identity}) yields
\[
-\frac{2}{n^2-1}y^i{\bf H}_{jl}+\frac{2}{n+1}y^i{\bf H}_{jl}=\frac{2(n-2)}{n^2-1}y^i{\bf H}_{jl}=0,
 \]
and f\/inally ${\bf H}_{jl}=0$, since $n\geq3$. Now, it results that $(M,F)$ is of constant f\/lag curvature, since ${\bf H}=0$. Conversely, suppose that
$(M,F)$ is of constant f\/lag curvature. Then, ${\bf H}=0$, $K_{kl}=K_{lk}$ and from (\ref{W H}) it follows that $\overset{*}{W}=W=0$, since $(M,F)$ is
of constant (scalar) f\/lag curvature.
\end{proof}

\begin{remark}
Projectively related Finsler metrics certainly have the same Weyl and Douglas curvatures. In \cite{ShenYu2008}, the authors studied projectively related Randers metrics. Their discussion is closely related to the subject of the present paper.
\end{remark}

\section{Projective vector f\/ields on Randers spaces}

Every vector f\/ield $V$ on $M$ induces naturally a transformation under the following inf\/initesimal coordinate transformations on $TM$,
$(x^i,y^i)\longrightarrow(\bar{x}^i,\bar{y}^i)$ given by
\[
\bar{x}^i=x^i+V^idt,\qquad \bar{y}^i=y^i+y^k\frac{\partial V^i}{\partial x^k}dt.
\]
This leads to the notion of \textit{the complete lift} $\hat{V}$ (or traditionally denoted by $V^C$, see \cite{Yano1957}) of $V$ to a vector f\/ield on $TM_0$ given by
\[
\hat{V}=V^i\frac{\partial}{\partial x^i}+y^k\frac{\partial V^i}{\partial x^k}\frac{\partial}{\partial y^i}.
\]
Since almost every geometric object in Finsler geometry depend on the both points and velocities, the Lie derivatives of such geometric objects should
be regarded with respect to $\hat{V}$. One may get familiar to the theory of Lie derivatives in Finsler geometry in \cite{Yano1957}. It is a notable remark
in the Lie derivative computations that ${\cal L}_{\hat{V}}y^i=0$ and the dif\/ferential operators ${\cal L}_{\hat{V}}$, $\frac{\partial}{\partial x^i}$
and~$\frac{\partial}{\partial y^i}$ commute.

A smooth vector f\/ield $V$ on $(M,F)$ is called \textit{projective} if each local f\/low dif\/feomorphism associated with $V$ maps geodesics onto geodesics. If $V$ is projective and each such map preserves af\/f\/ine parameters, then $V$ is called \textit{affine}, otherwise it is said to be \textit{proper
projective}. The collection of all projective vector f\/ields on $M$ is a f\/inite-dimensional Lie algebra, with respect to the usual Lie bracket operation on vector f\/ields, called the projective algebra, and is denoted by~$p(M,F)$. It is easy to prove that a vector f\/ield $V$ on the Finsler space $(M,F)$ is a projective if and only if there is a function $P$ def\/ined on $TM_0$ such that
\begin{gather}
\label{Lie}
{\cal L}_{\hat{V}}G^i=Py^i,
\end{gather}
and $V$ is af\/f\/ine if and only if $P=0$. Whence $F$ is Riemannian, the equation (\ref{Lie}) is just ${\cal L}_{\hat{V}}\Gamma^i_{jk}=\omega_j\delta^i_{\
k}+\omega_k\delta^i_{\ j}$, where $\omega_j$ are the components of a globally def\/ined 1-form on $M$ and thus, $P(x,y)=\omega_i(x)y^i$.

%\subsection{Proof of Theorem \ref{Liethm}}
\begin{proof}[Proof of Theorem \ref{Liethm}]
Suppose that $V$ is $F$-projective. Hence it preserves the Douglas tensor, i.e.\ ${\cal L}_{\hat{V}}D^i_{\ jkl}=0$. Let us put $T^i=\alpha s^i_{\
\circ}$. The sprays $G^i$ of $F$ and $\widehat{G}^i=G^i_\alpha+T^i$ are projectively related and thus they have the same Douglas tensor, hence
\[
D^i_{\ jkl}=\widehat{D}^i_{\ jkl}=\frac{\partial^3}{\partial y^j\partial y^k\partial y^l}\left\{T^i-\frac{1}{n+1}T^m_{\ y^m}y^i\right\}.
\]
A simple calculation shows that $T^m_{\ m}=0$. From that we have
\[
{\cal L}_{\hat{V}}D^i_{\ jkl}={\cal L}_{\hat{V}}T^i_{.j.k.l}={\cal L}_{\hat{V}}\{\alpha s^i_{\ \circ}\}_{.j.k.l}=0.
\]
Therefore, there are functions $H^i(x,y),\ (i=1,2,\dots,n)$ quadratic in $y$ such that
\begin{gather}
\label{Wee}
{\cal L}_{\hat{V}}\{\alpha s^i_{\ \circ}\}=H^i.
\end{gather}
Let us put $t_{ij}={\cal L}_{\hat{V}}a_{ij}$. Observe that ${\cal L}_{\hat{V}}\{\alpha s^i_{\ \circ}\}=\frac{t_{\circ\circ}}{2\alpha}s^i_{\
\circ}+\alpha{\cal L}_{\hat{V}} s^i_{\ \circ}$. Now the equation (\ref{Wee}) can be re-written as follows:
\begin{gather}
\label{Hii}
t_{\circ\circ}s^i_{\ \circ}+2\alpha^2{\cal L}_{\hat{V}}s^i_{\ \circ}=\alpha H^i.
\end{gather}
Here we emphasis that $\alpha^2=a_{ij}(x)y^iy^j$, $t_{\circ\circ}s^i_{\ \circ}=(t_{ij}(x)s^i_{\ k}(x))y^iy^ky^k$ and ${\cal L}_{\hat{V}}s^i_{\ \circ}=({\cal L}_{V}s^i_{\ k})(x)y^k$ are polynomials in $y^1,y^2,\dots,y^n$. Hence the left hand of (\ref{Hii}) is a polynomial in $y^1,y^2,\dots,y^n$ for every $i$, while the right hand is not. It follows immediately that $H^i=0$ for every index $i$ and~(\ref{Wee}) reads as ${\cal L}_{\hat{V}}\{\alpha s^i_{\ \circ}\}=0$. Recall that the geodesic coef\/f\/icients of $F$ are of the following form:
\begin{gather}
\label{Randers gco}
G^i = G^i_\alpha + \left({e_{\circ\circ} \over 2F} -s_\circ\right)y^i+ \alpha s^i_{\ \circ}.
\end{gather}
From ${\cal L}_{\hat{V}}\{\alpha s^i_{\ \circ}\}=0$ and ${\cal L}_{\hat{V}}G^i=Py^i$ it results now, that we have
\[
{\cal L}_{\hat{V}}G^i={\cal L}_{\hat{V}}\left\{G^i_\alpha + \left({e_{\circ\circ} \over 2F} -s_\circ\right)y^i\right\}=Py^i,
\]
and f\/inally we obtain
\[
{\cal L}_{\hat{V}}G^i_\alpha=\left\{P-{\cal L}_{\hat{V}}\left({e_{\circ\circ} \over 2F} -s_\circ\right)\right\}y^i,
\]
which shows that $V$ is a $\alpha$-projective vector f\/ield. Conversely suppose that $V$ is $\alpha$-projective (i.e.\
${\cal L}_{\hat{V}}G^i_\alpha=\omega_\circ y^i$, for some 1-forms $\omega_\circ=\omega_k(x)y^k$ on $M$)  and ${\cal L}_{\hat{V}}\{\alpha s^i_{\ \circ}\}=0$. From (\ref{Randers gco}) it follows
\begin{gather*}
{\cal L}_{\hat{V}}G^i = {\cal L}_{\hat{V}}\left\{G^i_\alpha+\left({e_{\circ\circ} \over 2F} -s_\circ\right)y^i+ \alpha s^i_{\ 0}\right\}={\cal L}_{\hat{V}}G^i_\alpha+{\cal L}_{\hat{V}}\left({e_{\circ\circ} \over 2F} -s_\circ\right)y^i \\
\phantom{{\cal L}_{\hat{V}}G^i}{}
 = \left\{\omega_\circ+{\cal L}_{\hat{V}}\left({e_{\circ\circ} \over 2F} -s_\circ\right)\right\}y^i,
\end{gather*}
which proves that $V$ is a $F$-projective vector f\/ield.
\end{proof}

%\subsection{Proof of Theorem \ref{admitting}}
Let us prove initially the following lemma:
\begin{lemma}
\label{lemkomaki}
Let $(M, F=\alpha+\beta)$ be an n-dimensional Randers space. If $s^i_{\ j}\neq0$, then $V$
is F-projective vector field if and only if it is a $\alpha$-homothety and ${\cal L}_{\hat{V}}d\beta = \mu d\beta$,
where ${\cal L}_{\hat{V}} a_{ij} = 2\mu a_{ij}$.
\end{lemma}

\begin{proof}
Suppose that $s^i_{\ \circ}\neq0$. By Theorem \ref{Liethm}, $V$ is $F$-projective if and only if it is $\alpha$-projective and
${\cal L}_{\hat{V}}\{\alpha s^i_{\ \circ}\}=0$. Let us suppose $t_{ij}={\cal L}_{\hat{V}}a_{ij}$ and observe that ${\cal L}_{\hat{V}}\{\alpha s^i_{\
\circ}\}=0$ is equivalent to
\begin{gather}
\label{tt}
t_{\circ\circ}s^i_{\ \circ}+2\alpha^2{\cal L}_{\hat{V}}s^i_{\ \circ}=0.
\end{gather}
It follows that $\alpha^2$ divides $t_{\circ\circ}s^i_{\ \circ}$ for every index $i$. This equivalent to that $s^i_{\ j}=0$ or $\alpha^2$ divides
$t_{\circ\circ}$ which means that $V$ is a conformal vector f\/ield on $(M,\alpha)$, since $s^i_{\ j}\neq0$. Since $V$ is already $\alpha$-projective,
it follows that $V$ is an $\alpha$-homothety and there is a constant $\mu$ such that ${\cal L}_{V}a_{ij}=2\mu a_{ij}$. From (\ref{tt}) we obtain
${\cal L}_{V}s^i_{\ j}=-\mu s^i_{\ j}$. Now observe that
\begin{gather*}
{\cal L}_{V}s_{ij} = {\cal L}_{V}\{a_{ik}s^k_{\ j}\}=({\cal L}_{V}a_{ik})s^k_{\ j}+a_{ik}{\cal L}_{V}s^k_{\ j}=2\mu a_{ik}s^k_{\ j}-\mu a_{ik}s^k_{\ j}=\mu s_{ij}.\tag*{\qed}
\end{gather*}
\renewcommand{\qed}{}
\end{proof}

\begin{proof}[Proof of Theorem \ref{admitting}] Let us suppose that $M$ is compact and $F=\alpha+\beta$ is a Randers metric of constant f\/lag curvature and $n\geq3$. Following \cite{BaoRobles2003}, $F$ is of constant {\bf S}-curvature and due to a result about constant {\bf S}-curvature Finsler spaces in \cite{Mo2008}, it follows that ${\bf S}=0$. This results $e_{\circ\circ}=r_{\circ\circ}+2\beta s_\circ=0$. Now let us suppose that $s^i_{\ j}\neq0$. By Lemma~\ref{lemkomaki}, every $F$-projective vector f\/ield $V$ is an $\alpha$-homothety and since $M$ is compact, thus every $F$-projective vector f\/ield~$V$ is $\alpha$-Killing. Hence in this case we have the inclusion $p(M,F)\subseteq k(M,\alpha)$, where $k(M,\alpha)$ denotes the Lie algebra of $\alpha$-Killing vector f\/ields. It is well-known that the dimension of algebra of $\alpha$-Killing vector f\/ields is at most $\frac{n(n+1)}{2}$. Therefore  $\dim(p(M,F))\leq\frac{n(n+1)}{2}$. Now let us suppose that $s^i_{\ j}=0$. In this case we have $p(M,F)=p(M,\alpha)$ and moreover one conclude that $\nabla_jb_i=0$ and~$F$ is a Berwald metric. Since $F$ is of constant f\/lag curvature, thus~$F$ and~$\alpha$ are metrics of zero f\/lag curvatures. Notice that $F$ is of constant f\/lag curvature, its Weyl curvature vanishes and since $s^i_{\ j}=0$, thus $F$ and $\alpha$ are projectively related and hence $\alpha$ has vanishing Weyl curvature and by Beltrami's theorem, $\alpha$ has constant sectional curvature. It is well-known that the dimension of $p(M,\alpha)$ is $n(n+2)$. Hence  we have $\dim  (p(M,F) )=n(n+2)$.
\end{proof}

The following inclusive result follows from the proof of Theorem \ref{admitting}.
\begin{corollary}
Let $(M,F=\alpha+\beta)$ be a Randers space of constant flag curvature. The following statements hold:
\begin{enumerate}\itemsep=0pt
\item[$(a)$] if $\beta$ is closed, then $p(M,F)=p(M,\alpha)$;
\item[$(b)$] if $\beta$ is not a closed $1$-form, then $p(M,F)\subseteq h(M,\alpha)$, where $h(M,\alpha)$ denotes the Lie algebra of $\alpha$-homothety vector fields.
    \end{enumerate}
\end{corollary}

%\subsection{Proof of Theorem \ref{mainthm}}
\begin{proof}[Proof of Theorem \ref{mainthm}]
To obtain general formulae, let us assume that $(M,F=\alpha+\beta)$ be a~Randers space of isotropic ${\bf S}$-curvature ${\bf S}=(n+1)\sigma(x)F$ and $V$ be a~non-af\/f\/ine projective vector f\/ield. Following a result in \cite{ShenXing2008}, it results that $e_{\circ\circ}=2\sigma(x)(\alpha^2-\beta^2)$. Suppose that there is a~function $\Psi$ such that $\Psi(x,y)$ is linear with respect to $y$ such that ${\cal L}_{\hat{V}}G^i=\Psi y^i$. By applying Theorem \ref{Liethm}, we have
\[
{\cal L}_{\hat{V}}G^i={\cal L}_{\hat{V}}\widetilde{G}^i+{\cal L}_{\hat{V}}(\sigma(\alpha-\beta)y^i)-{\cal L}_{\hat{V}}s_\circ y^i=\Psi y^i.
\]
Put $t_{ij}={\cal L}_{\hat{V}}a_{ij}$. It is well-known that ${\cal L}_{\hat{V}}y^i=0$, it follows $t_{\circ\circ}={\cal L}_{\hat{V}}\alpha^2$ and
\[
{\cal L}_{\hat{V}}\widetilde{G}^i+\alpha{\cal L}_{\hat{V}}\sigma y^i-\beta{\cal L}_{\hat{V}}\sigma y^i
+\frac{t_{\circ\circ}}{2\alpha}cy^i-\sigma y^i{\cal L}_{\hat{V}}\beta
-{\cal L}_{\hat{V}}s_\circ y^i =\Psi y^i.
\]
Given any natural coordinate system $((x^i,y^i),\pi^{-1}(U))$ and $x\in U$, we can regard each terms of the above equation as a polynomial in
$y^1,y^2,\dots,y^n$. Multiplying the two sides of the last equation in $\alpha$, we obtain the following relation:
\[
{\rm Rat}^i+\alpha\, {\rm Irrat}^i=0,\qquad i=1,2,\dots,n,
\]
where the polynomials ${\rm Rat}^i$ and ${\rm Irrat}^i$ are given by
\begin{gather*}
%\label{Rat}
{\rm Rat}^i=\alpha^2{\cal L}_{\hat{V}}\sigma y^i+\frac{1}{2}t_{\circ\circ}\sigma y^i,\\
%\label{Irrat}
{\rm Irrat}^i={\cal L}_{\hat{V}}\widetilde{G}^i-(\beta{\cal L}_{\hat{V}}\sigma+\sigma{\cal L}_{\hat{V}}\beta+{\cal L}_{\hat{V}}s_\circ +\Psi) y^i.
\end{gather*}
Now let us assume ${\bf S}=(n+1)\sigma(x)F=0$. By Lemma \ref{lemkomaki}, $(M,F)$ must be locally projectively f\/lat, otherwise $V$ is a $\alpha$-homothety
which is a contradiction to the assumption that $V$ is non-$\alpha$-homothety. Hence $s_{ij}=0$ and by $e_{ij}=r_{ij}+b_is_j+b_js_i=r_{ij}=0$. This is
equivalent to $\nabla_ib_j=0$ and $(M,F)$ is a Berwald space.
\end{proof}

By Theorem \ref{admitting}, the following theorems results
\begin{theorem}
\label{fini}
Let $(M,F=\alpha+\beta)$ be a compact $n$-Randers space of constant flag curvature. The following statements hold:
\begin{enumerate}\itemsep=0pt
\item[$(a)$] if $\dim (p(M,F) )=\frac{n(n+1)}{2}$, then $\alpha$ is of constant sectional curvature;
\item[$(b)$] if $\dim (p(M,F) )>\frac{n(n+1)}{2}$, then $F$ is a locally Minkowski metric.
\end{enumerate}
\end{theorem}

\begin{proof}
Let us suppose $\dim(p(M,F))=\frac{n(n+1)}{2}$. Due to the discussions in the proof of Theorem~\ref{admitting}, in this case we have
$p(M,F)\subseteq k(M,\alpha)$ and thus,
\[
\frac{n(n+1)}{2}=\dim (p(M,F))\leq \dim(k(M,\alpha)).
\]
Hence $(M,\alpha)$ is of maximum rank $\frac{n(n+1)}{2}$ and it is well-known that in this case $\alpha$ is of constant sectional curvature. This
proves $(a)$. To prove $(b)$, we notice that if $\dim  (p(M,F) )>\frac{n(n+1)}{2}$, then we must have $s^i_{\ j}=0$, Otherwise, by proof of
Theorem \ref{admitting}, we have $\dim (p(M,F) )\leq\frac{n(n+1)}{2}$. Now, $s^i_{\ j}=0$ and ${\bf S}=0$ results that $F$ is a Berwald
metric which is already of constant f\/lag curvature. $F$ is not Riemannian and Numata's theorem ensures that ${\bf K}=0$. Finally, Akbar-Zadeh's
classif\/ication theorem entails $F$ is a locally Minkowski metric.
\end{proof}

\begin{example}
In \cite{BaoShen2002}, the authors presented a worthily source of a 1-parameter family of Randers metric $F_\kappa=\alpha_\kappa+\beta_\kappa$ on the Lie group $S^3$ which all are of constant positive f\/lag curvature~$\kappa$. Due to their construction, non of the Riemannian metrics $\alpha_\kappa$ is of constant sectional curvature and hence, by Theorem \ref{fini}, it follows that $\dim\big(p(S^3,F_\kappa)\big)<6$.
\end{example}

\section{The Lorentz nonlinear connection\\ and Randers projective symmetry}

The stage on which special relativity is played out is a specif\/ic four dimensional manifold, known as
Minkowski spacetime. The $x^\mu$, $\mu=0,1,2,3$, are coordinates on this manifold and conventionally, we set $x^0=t$. The elements of spacetime are
known as events; an event is specif\/ied by giving its
location in both space and time. The inf\/initesimal (distance) between two points known as the \textit{spacetime interval} is def\/ined by
\[
ds^2=\eta_{\mu\nu}dx^\mu dx^\nu=-dt^2+dx^2+dy^2+dz^2.
\]
The matrices $\Lambda$ satisfying $\Lambda^T\eta\Lambda=\eta$ are known as the \textit{Lorentz transformations}. As a notable well-known case, consider the celebrated Randers metric of the form $F=\sqrt{\eta_{\mu\nu}dx^\mu dx^\nu}+{\bf A}_idx^i$ on the 4-manifold
of spacetime, where ${\bf A}$ is the electromagnetic vectorial potential and ${\bf F}=d{\bf A}$ obtained in the Cartesian coordinates $(t,x,y,z)$ as
\begin{gather*}
{\bf F}_{\mu\nu}=
  \begin{pmatrix}
    0 & -E_x & -E_y & E_z \\
    E_x & 0 & B_z & -B_y \\
    E_y & -B_z & 0 & B_x \\
    E_z & B_y & -B_x & 0
  \end{pmatrix}.
\end{gather*}

Notice that ${\bf F}=0$ if and only if $F$ is locally projectively f\/lat. Hence, the presence of a proper electromagnetic f\/ield, entails non-locally
projectively f\/latness of $F$. Consider a Lorentz transformation $\Lambda$ which maps the coordinates $(t,x,y,z)$ onto
$(\bar{t},\bar{x},\bar{y},\bar{z})$. $\Lambda$ changes $F=\sqrt{\eta_{\mu\nu}dx^\mu dx^\nu}+{\bf A}_idx^i$ to $\bar{F}=\sqrt{\eta_{\mu\nu}dx^\mu
dx^\nu}+\bar{{\bf A}}_idx^i$. Following an extensive discussion on projectively related Randers metrics in~\cite{ShenYu2008}, we conclude that the
Lorentz transformation~$\Lambda$ is $F$-projective if and only if $\Lambda^T{\bf F}\Lambda={\bf F}$. The collection of all such Lorentz transformation
forms a subgroup of the Lorentz group which is at once a subgroup of projective group $P(M,F)$.
%*****************************************************************************************
%\subsection{Lorentz nonlinear connection and Randers projective symmetry}
Theory of Finsler spaces with $(\alpha,\beta)$-metrics was studied by such famous geometers as M.~Matsumoto, D.~Bao and many others as a natural
extension of the theory of Randers spaces. Associated to any $(\alpha,\beta)$-metric $F=F(\alpha,\beta)$ one may consider a nonlinear connection called
Lorentz connection which has physical applications in the study of gravitational and electromagnetic f\/ields~\cite{Miron2006}. In this section, we uncover
some results about its projective symmetry in a Finsler space with a~Randers metric.

Let $F=F(\alpha,\beta)$ be an $(\alpha,\beta)$-metric on the manifold $M$. Through a variational approach, Lorentz equations are derived using
Euler--Lagrange equations in the following form:
\[
\frac{d^2x^i}{ds^2}+\Gamma^i_{jk}\frac{dx^j}{ds}\frac{dx^k}{ds}+\sigma\left(x,\frac{dx}{ds}\right)s^i_{\ j}\frac{dx^j}{ds}=0,
\]
where $ds^2=\alpha^2(x,\frac{dx}{dt})dt^2$ and $\sigma=\alpha{F^2_\beta}/{F^2_\alpha}$. The \textit{Lorentz nonlinear connection}
$\overset{\circ}{G}{^i_j}$ is now def\/ined~by
\begin{gather*}
%\label{Lorentz connection}
\overset{\circ}{G}{^i_{\ j}}(x,y)=\Gamma^i_{jk}(x)y^k+\sigma(x,y)s^i_{\ j}.
\end{gather*}
Every geometric object associated to the Lorentz connection will be denoted by the sign ``$^\circ$'' on top. Notice that the Lorentz nonlinear
connection is determined by the Finsler--Lagrange metric $F=F(\alpha,\beta)$ only. Notice that the autoparallel curves of the nonlinear connection
$\overset{\circ}{G}{^i_j}$, in the natural parameterizations (i.e.~$\alpha(x, \frac{dx}{ds}) = 1$), coincide with the integral curves of the canonical
semispray $S$ given by
\begin{gather*}
%\label{Lorentz spray}
S=y^i\frac{\partial}{\partial x^i}-2\overset{\circ}{G}{^i}\frac{\partial}{\partial y^i},\qquad \textrm{where} \qquad
2\overset{\circ}{G}{^i}=\Gamma^i_{jk}y^jy^k+\alpha s^i_{\ \circ}.
\end{gather*}

Akbar-Zadeh in \cite{Akbar-Zadeh1995} considers the Berwald connection of the semispray $\overset{\circ}{G}{^i}$ to obtain a covariant derivative and derived a unif\/ied formulation for electromagnetism and gravity. However, it encounters physical consistency: all the formulation require being invariant under the Lorentz group in f\/lat space-time. This is not satisf\/ied generally. It can be shown that the only Lorentz transformation which preserve the spray $\overset{\circ}{G}{^i}$ are those satisfying~$\Lambda^T{\bf F}\Lambda={\bf F}$.

\subsection*{Acknowledgements}

The authors would like to express their grateful thanks to the referees for their valuable comments.
This work was supported in part by the Institute for Research in Fundamental
Sciences (IPM) by the grant No.~89530036. The f\/irst author wishes to thank
Borzoo Nazari for many fruitful conversations.

\vspace{-2mm}

\pdfbookmark[1]{References}{ref}
\LastPageEnding

\end{document}